\documentclass[12pt,twoside]{amsart}

\usepackage{amsmath}
\usepackage{amssymb}
\usepackage{graphicx,epsf,amsmath}  
\usepackage{epsf,graphicx}
\usepackage{comment}

\newtheorem{thm}{Theorem}[section]
\newtheorem{lemma}[thm]{Lemma}

\newtheorem*{thm*}{Theorem}

\newcommand{\lp}{\left(}
\newcommand{\rp}{\right)}

\newcommand{\LL}{{\mathcal L}}
\newcommand{\EE}{{\mathcal E}}

\newcommand{\R}{\mathbb{R}}

\newcommand{\abs}[1]{\left\lvert #1\right\rvert}

\newcommand{\be}{\begin{equation}}
\newcommand{\ee}{\end{equation}}
\newcommand{\bee}{\begin{equation*}}
\newcommand{\eee}{\end{equation*}}
\newcommand{\bea}{\begin{eqnarray}}
\newcommand{\eea}{\end{eqnarray}}
\newcommand{\bs}{\begin{split}}
\newcommand{\es}{\end{split}}
\renewcommand{\le}         {\leqslant}
\renewcommand{\leq}         {\leqslant}
\renewcommand{\ge}         {\geqslant}
\renewcommand{\geq}         {\geqslant}
\renewcommand{\div}         {\,{\rm{div}}}
\title[Elliptic PDEs on manifolds with boundary ]{Some 
elliptic PDEs\\
on Riemannian manifolds with boundary}

\author{Yannick Sire}
\author{Enrico Valdinoci}\thanks{
{\it YS}:
Universit\'e Aix-Marseille 3, Paul C\'ezanne --
LATP --
Marseille, France 
and Laboratoire Poncelet,
UMI 2615--
Moscow, Russia--
{\tt sire@cmi.univ-mrs.fr}\\
{\it EV}:
Universit\`a di Roma Tor Vergata --
Dipartimento di Matematica --
Rome, Italy --
{\tt enrico.valdinoci@uniroma2.it}
}
\begin{document}
\begin{abstract}
The goal of this paper is to investigate some rigidity properties 
of stable solutions of elliptic equations set on manifolds with 
boundary. 

We provide several types of 
results, according to the dimension of
the manifold and the sign of its Ricci curvature.   
\end{abstract}
\maketitle
\tableofcontents

\section{Introduction}

Let $(\mathcal M,\bar g)$ be a complete,
connected,
smooth, $n+1$-dimensional
manifold with boundary $\partial \mathcal M$,
endowed with a smooth Riemannian metric $\bar g=\{ \bar g_{ij}\}_{i,j=1,...,n}$.

The volume element writes in local coordinates as
\begin{equation}\label{vg}
dV_{\bar g} = \sqrt{|\bar g|}\,dx^1\wedge \dots \wedge dx^n,
\end{equation}
where $\{ dx^1,\dots,dx^n\}$ is the basis of $1$-forms dual to the vector basis $\{\partial_i,\dots,\partial_n\}$ and we use the standard notation $|\bar g|=\det(\bar g_{ij})\ge0$.

We denote by $\div_{\bar g} X$ the divergence of a smooth vector
field $X$ on $\mathcal M$, that is, in local coordinates,
\begin{equation*}
\div_{\bar g} X=\frac{1}{\sqrt{|\bar g|}} 
\partial_i \Big( \sqrt{|\bar g|}X^i\Big),\end{equation*}
with the Einstein summation convention.

We also denote by $\nabla_{\bar g}$ the Riemannian gradient
and by $\Delta_{\bar g}$ the Laplace-Beltrami operator, that is,
in local coordinates,
\begin{equation}
\label{1.0}
(\nabla_{\bar g}\phi)^i=\bar g^{ij}\partial_j \phi
\end{equation}
and
\begin{equation*}
\Delta_{\bar g} \phi= \div_{\bar g}(\nabla_{\bar g} \phi)
=
\frac{1}{\sqrt{|\bar g|}}\partial_i \Big(
{\sqrt{|\bar g|}} \bar g^{ij}\partial_j \phi\Big)
,\end{equation*}
for any smooth function $\phi:\mathcal M\rightarrow \R$.


We set~$\langle \cdot,\cdot\rangle$ to be the scalar product
induced by $\bar g$. 

Given a vector field $X$, we also denote
$$ |X|=\sqrt{\langle X,X\rangle}.$$

Also (see, for instance Definition~3.3.5 in~\cite{Jost}),
it is customary to define
the Hessian of a smooth function $\phi$ as
the symmetric $2$-tensor given in a local patch by
$$ (H_{\bar g} \phi)_{ij}=\partial^2_{ij}\phi-\Gamma^k_{ij}\partial_k\phi,$$
where $\Gamma^k_{ij}$ are the  
Christoffel symbols, namely
$$ 
\Gamma_{ij}^k=\frac12 \bar g^{hk} \left( \partial_i
\bar g_{hj} +\partial_j 
\bar g_{ih} -\partial_h \bar g_{ij} \right).$$
Given a tensor $A$,
we define its norm by $|A|=\sqrt{A A^*}$, where $A^*$
is the adjoint.

The present
paper is devoted to the study of special solutions of elliptic 
equations on manifolds with boundary and is, in some sense, a follow up to the paper by
the authors and Farina (see \cite{FSV-man}) where the case without 
boundary was investigated. In an Euclidean context,
i.e. $\mathcal M=\R^{n+1}_+$ with the flat metric, the
rigidity features of the 
stable solutions has been investigated in \cite{SV1,CS}.

Boundary problems are related (via a theorem of Caffarelli and 
Silvestre \cite{cafS}) to non local equations involving fractional 
powers of the Laplacian. An analogue of the results 
of~\cite{cafS}
has been obtained in a geometric context, by means of scattering theory (see \cite{FG,GJMS,GZ}). 

In this paper, we will focus on the following two specific models: 
\begin{itemize}
\item {\bf product manifolds} of the type
$$\Big ( \mathcal M= M \times \R^+, \bar g=g +|dx|^2 \Big )$$
where $(M,g)$ is a complete, smooth
Riemannian manifold
without boundary,
and \item {\bf the hyperbolic halfspace}, i.e.
$$\Big ( \mathcal M= \mathbb H^{n+1}, 
\bar g=\frac{|dy|^2+|dx|^2}{x^2} \Big )$$
where $x>0$ and $y \in \R^n$. 
\end{itemize}
Notice that the above models comprise both the
positive and the negative 
curvature cases.

We denote by~$\nu$ the exterior derivative at points
of~$\partial{\mathcal M}$.

We will investigate the two 
following problems
\begin{equation}\label{pbpos}
\left \{
\begin{array}{c}
 \Delta_{\bar g} u=0\,\,\,\mbox{ in }\, 
\mathcal M=
M \times \R^+,\\
\partial_\nu u=f(u)\,\,\,\,\mbox{ on }\,M \times  \left \{0 \right 
\}. 
\end{array} \right . 
\end{equation}
and 
\begin{equation}\label{pbneg}
\left \{
\begin{array}{c}
-\Delta_{\bar g} u-s(n-s)u=0\,\,\,\mbox{ in }\, 
\mathcal M=
\mathbb H^{n+1},\\
\partial_\nu u=f(u)\,\,\,\,\mbox{ on }\, \partial \mathbb H^{n+1}. 
\end{array} \right . 
\end{equation}
where $f$ is a $C^1(\mathcal M)$ nonlinearity
(in fact, up to minor modifications, the proofs
we present also work for
locally Lipschitz nonlinearities).

The real parameter $s$ in~\eqref{pbneg}
is chosen to be
$$s=\frac{n}{2}+\gamma,$$ where 
$\gamma \in (0,1)$. 

We recall that the problem in~\eqref{pbpos}
has been
studied in the context of conformal geometry 
and it is related to
conformally compact Einstein manifolds
(see section~\ref{Ein} below
for a further 
discussion).

We will consider weak solutions of \eqref{pbpos} and 
\eqref{pbneg}.
Namely, we say that $u$ is weak solution of \eqref{pbpos} if,
for every $\xi \in C_0^\infty(M\times \R)$, we have
\begin{equation}\label{weakformpos}
\int_{\mathcal M } \langle \nabla_{\bar g} u, \nabla_{\bar g} \xi \rangle=\int_{\partial \mathcal M} f(u) \xi. 
\end{equation}
Analogously, we say that $u$ is weak solution of 
\eqref{pbneg} 
if, for every $\xi \in C^\infty(\mathcal M)$, we have
that 
\begin{equation}\label{weakformneg}
\int_{\mathcal M } \langle \nabla_{g} u, \nabla_{g} \xi \rangle - s(n-s)\int_{\mathcal M } u \xi =\int_{\partial \mathcal M} f(u) \xi. 
\end{equation}

We focus on an important class of solutions of \eqref{pbpos} and 
\eqref{pbneg}, namely the so called stable solutions. 

These 
solutions play an 
important role in the calculus of variations and are characterized by 
the fact that the second variation of the energy functional
is non negative definite.
This condition may be explicitly written in our case by saying
that a solution~$u$ of either~\eqref{pbpos} and
or \eqref{pbneg} is stable if
\begin{equation}\label{STA}
\int_{\mathcal M} |\nabla_{\bar g}\xi|^2 dV_{\bar g}-s(n-s)\varepsilon \int_{\mathcal M} \xi^2 -\int_{\partial \mathcal M} f'(u) \xi^2\,dV_{\bar g}
\ge0\end{equation}
for every~$\xi\in C^\infty_0 (M \times \R)$
with~$\varepsilon =0 $ 
in case of \eqref{pbpos},
and
for every~$\xi\in C^\infty( \mathcal M)$
and with~$\varepsilon =1$
in case of \eqref{pbneg}.

\subsection{Results for product manifolds}

Now we present our results in the case of
product manifolds~$\mathcal M= M \times \R^+$.

\begin{thm}\label{liouville1}
Assume that the metric on $\mathcal M= M \times \R^+$ is given by
$\bar g= g+ |dx|^2$. 

Assume furthermore that $M$ is compact and satisfies 
$$Ric_g \geq 0$$
with~$Ric_g$
not vanishing identically. 

Then every bounded stable weak solution~$u$ of~\eqref{pbpos} is 
constant. 
\end{thm} 

We remark that the assumption on the boundedness of~$u$ is needed 
as the following example shows: the function $u(x,y)=x$
is a stable solution of 

\begin{equation*}
\left \{
\begin{array}{c}
\Delta_{\bar g} u=0\,\,\,\mbox{in}\, M \times \R^+,\\
\partial_\nu u=-1\,\,\,\,\mbox{on}\,M \times  \left \{0 \right \}. 
\end{array} \right. 
\end{equation*}

{F}rom theorem~\ref{liouville1}, one also obtains
the following Liouville-type theorem for the 
half-Laplacian on compact manifolds (for the definition
and basic functional properties of fractional
operators see, e.g.,~\cite{Kato}):

\begin{thm}\label{liouville1NL}
Let $(M,g)$ be a compact manifold and
$u:M  \to \R$ be a smooth bounded solution of 
\begin{equation}\label{uhdsiuuiuu111}
(-\Delta_g)^{1/2}u=f(u),\end{equation}
with
\begin{equation}\label{staze}
\int_{ \mathcal M }  ( |\nabla_g \xi|^2+|\nabla_x \xi|^2 )
- \int_{\partial\mathcal M} f'(u)\xi^2 \ge 0,\end{equation}
for every~$\xi\in C^\infty_0 ({\mathcal{M}})$.

Assume furthermore that 
$$Ric_g \geq 0$$
and~$Ric_g$ does not vanish identically. 

Then $u$ is constant. 
\end{thm} 

Results for~$(-\Delta_g)^\alpha$ with~$\alpha\in(0,1)$
may be obtained with similar techniques as well.

\begin{thm}\label{liouville2}
Assume that the metric on $\mathcal M
= M \times \R^+$ is given by $\bar g= g+ 
|dx|^2$,
that $M$ is complete, and
$$Ric_g \geq 0,$$
with~$Ric_g$ not vanishing identically. 

Assume also that, for any $R>0$, the volume of
the geodesic ball~$B_R$ in~$M$ (measured with respect
to the volume element $dV_g$) 
is bounded by~$C( R+1)$, for some~$C>0$.

Then every bounded stable weak solution
$u$ of \eqref{pbpos} is constant. 
\end{thm} 

Next theorem is a flatness result when the Ricci tensor of $M$ 
vanishes identically:

\begin{thm}\label{flatness}
Assume that the metric on~$\mathcal M= M \times \R^+
$ is given by $\bar g= 
g+ 
|dx|^2$ and~$Ric_g$ vanishes identically.

Assume also that, for any $R>0$, the volume of
the geodesic ball~$B_R$ in~$M$ (measured with respect
to the volume element $dV_g$)
is bounded by~$C( R+1)$, for some~$C>0$.

Then for every $x>0$ and~$c\in\R$,
every connected component of the 
submanifold
$$\mathcal S_x = \{ y\in M,\,\,\,u(x,y)=c\}$$
is a geodesic. 
\end{thm} 

\subsection{Results for the hyperbolic space}

The next theorem provides a flatness result when the manifold $\mathcal M$
is $ \mathbb H^{3}$.

\begin{thm}\label{thhyp}
Let~$n=2$.

Let $u$ be a smooth weak solution of~(\ref{pbneg}) and let 
$s=\frac{n}{2}+\gamma$ where $\gamma \in (0,1)$.


Also, suppose that either
\begin{equation}\label{001} \partial_{y_2} u>0
\end{equation}
or 
\begin{equation}\label{002}f' \leq 
0\,\,\,\,\mbox{on}\,\,\,\partial \mathbb H^{n+1}.
\end{equation}

Then, for every $x>0$ and~$c\in\R$,
each of the submanifold  
$$\mathcal S_x =\left \{ y\in \R^n,\,\,\,|\,\,u(x,y)=c x^{n-s} 
\right \} $$
is a Euclidean straight line.
\end{thm}

As discussed in details in section~\ref{989},
the proof of theorem~\ref{thhyp}
contains two main ingredients: 
\begin{enumerate}
\item We first notice that the metric on $\mathbb H^{n+1}$ is conformal to the flat metric on $\R^{n+1}_+$. 
\item We then use some results by the authors in \cite{SV1} (see also
\cite{CS} for related problems) to 
get the desired result.  
\end{enumerate}

The rest of this paper is structured as 
follows. In section~\ref{S-1}
we prove a geometric inequality for stable
solutions in product manifolds, from which
we obtain the proofs of 
theorems~\ref{liouville1}--\ref{flatness}, contained in
section~\ref{s-2-2}. Then, in section~\ref{negative},
we consider the hyperbolic case and we prove 
theorem~\ref{liouville1}.

\section{The case of product manifolds and a
weighted Poincar\'e 
inequality for stable solutions of \eqref{pbpos}}\label{S-1}

Now we deal with the case of product manifolds~$
M \times \R^+$.

In order to simplify notations,
we write $\nabla $ instead
of $\nabla_{\bar g}$ for the gradient on
$M \times \R^+$ but we will keep
the notation $\nabla_g$ for the Riemannian gradient on $M$.

Recalling~\eqref{1.0}, we have that
\begin{equation}\label{0.1}
\nabla=(\nabla_g,\partial_x).\end{equation}
In the subsequent theorem~\ref{0090},
we obtain a formula involving the geometry, in a quite implicit 
way, of the 
level sets of stable solutions of \eqref{pbpos}. 

Such a formula may be considered
a geometric version of
the Poincar\'e
inequality, since the $L^2$-norm of the gradient
of any test function bounds the~$L^2$-norm of the test
function itself. Remarkably, these $L^2$-norms are
weighted and the weights have a neat geometric
interpretation.

These type of geometric
Poincar\'e
inequalities were first obtained 
by~\cite{SZ1,SZ2}
in the Euclidean setting, and similar estimates
have been recently widely used
for rigidity results in PDEs (see, for 
instance,~\cite{FScV,SV1,FER}).
 
\begin{thm}\label{0090}
Let $u$ be a stable solution of \eqref{pbpos}
such that $\nabla_g u$ is bounded.

Then, for every $\varphi \in C_0^{\infty}(M\times \R)$, the
following inequality holds:
\begin{equation}\label{poincare}
\begin{split}&
\int_{M \times \R^+} \Big \{ Ric_g (\nabla_g u, \nabla_g  u) +|H_g 
u|^2 -|\nabla_g |\nabla_g u||^2 \Big \} \varphi^2\\&\qquad\quad 
\leq 
\int_{M \times \R^+} |\nabla_g u |^2 | \nabla \varphi |^2.    
\end{split}\end{equation} 
\end{thm}

Notice that only the geometry of $M$ comes into play
in formula \eqref{poincare}.

\begin{proof}
First of all,
we recall
the classical Bochner-Weitzenb\"ock formula
for a smooth function~$\phi:{\mathcal{M}}\rightarrow\R$
(see,
for instance,~\cite{Berger, Wang} and references therein):
\begin{equation}\label{BOC}
\frac 12\Delta_{\bar g} |\nabla_{\bar g} \phi|^2=
|H_{\bar g} \phi|^2+ \langle \nabla_{\bar g} \Delta_{\bar g}
\phi,\nabla_{\bar g}\phi\rangle
+Ric_{\bar g} (\nabla_{\bar g}
\phi,\nabla_{\bar g}\phi).\end{equation}

The proof of theorem~\ref{0090}
consists in plugging the
test function~$\xi= |\nabla_g u| 
\varphi$ in the stability condition \eqref{STA}: after
a simple computation, this gives
\begin{equation}\label{AQ}
\begin{split} &
\int_{\mathcal M} \varphi^2   |\nabla |\nabla_g u | |^2 +\frac{1}{2} \langle \nabla |\nabla_g u|^2, \nabla \varphi^2 \rangle + |\nabla_g u|^2  
|\nabla \varphi |^2 \\ &\qquad-
\int_M f'(u)  |\nabla_g u|^2  \varphi^2 \geq 0. \end{split}
\end{equation}
Also, by recalling~\eqref{0.1}, we have
\begin{equation}\label{901}
\langle \nabla |\nabla_g u|^2, \nabla 
\varphi^2 \rangle= 
\langle \nabla_g |\nabla_g u|^2, \nabla_g \varphi^2 \rangle +
\partial_x |\nabla_g u|^2 \partial_x \varphi^2 .
\end{equation}
Moreover, since~$M$ is boundaryless,
we can use on~$M$ the Green formula
(see,
for example,
page~184 of~\cite{Gallot}) and obtain that
\begin{equation}\label{902}
\begin{split}
& \int_{\mathcal{M}} \langle \nabla_g |\nabla_g u|^2, \nabla_g 
\varphi^2 \rangle = \int_{\R^+} \int_M
\langle \nabla_g |\nabla_g u|^2, \nabla_g \varphi^2 \rangle
\\ &\qquad\quad= -\int_{\R^+}\int_{M}
\Delta_g |\nabla_g u| \varphi^2 =
-\int_{\mathcal M}
\Delta_g |\nabla_g u| \varphi^2.
\end{split}\end{equation}
Hence,
using~\eqref{BOC}, \eqref{901} and~\eqref{902},
we conclude that
\begin{equation}\begin{split}\label{AX}
& \frac{1}{2} \int_{\mathcal 
M}  \langle 
\nabla |\nabla_g u|^2, \nabla \varphi^2 \rangle= \frac{1}{2} 
\int_{\mathcal M} \partial_x |\nabla_g u|^2 \partial_x \varphi^2 
\\
&\qquad-
\int_{\mathcal M}  \varphi^2 \Big\{ |H_g u|^2 + \langle \nabla_g 
\Delta_g u , \nabla_g u \rangle +Ric_g (\nabla_g u , \nabla_g u) 
\Big\} . \end{split}\end{equation}
Using the equation in~\eqref{pbpos},
we obtain that
$$\Delta_g u = 
-\partial_{xx}u, $$
so~\eqref{AX} becomes
\begin{equation}\label{AY}
\begin{split}&
\frac{1}{2} \int_{\mathcal M}  \langle \nabla |\nabla_g u|^2, 
\nabla \varphi^2 \rangle= \frac{1}{2} \int_{\mathcal M} \partial_x 
|\nabla_g u|^2 \partial_x \varphi^2 \\&\qquad-
\int_{\mathcal M}  \varphi^2 |H_g u|^2 + 
\int_{\mathcal M} \varphi^2\langle \nabla_g \partial_{xx} u , 
\nabla_g u \rangle -  \int_{\mathcal M} \varphi^2 Ric_g (\nabla_g 
u , \nabla_g u)  . \end{split}\end{equation}
Furthermore,
Integrating by parts, we see that
\begin{eqnarray*}
&& \int_{\mathcal M}\partial_x |\nabla_g 
u|\partial_x\varphi^2=
\int_{M}\int_0^{+\infty}
\partial_x |\nabla_g 
u|\partial_x\varphi^2\\ &&\qquad=
-\int_M \Big( \partial_x |\nabla_g
u| \varphi^2\Big)|_{x=0}-\int_M 
\int_0^{+\infty}
\partial_{xx} |\nabla_g
u| \varphi^2\\
&&\qquad=
-\int_M \Big( \partial_x |\nabla_g
u| \varphi^2\Big)|_{x=0}-\int_{\mathcal{M}}
\partial_{xx} |\nabla_g
u| \varphi^2.
\end{eqnarray*}
Consequently, \eqref{AY} becomes
\begin{eqnarray}\label{AZ}\nonumber&& 
\frac{1}{2} \int_{\mathcal M}  \langle \nabla |\nabla_g u|^2, 
\nabla \varphi^2 \rangle=\\&&\qquad - \int_{\mathcal M} \varphi^2  
\Big \{ 
\frac{1}{2}\partial_{xx} |\nabla_g u|^2 + |H_g u|^2 +Ric_g 
(\nabla_g u , \nabla_g u) \Big \} \\ &&\qquad+
\int_{\mathcal M} \varphi^2  \langle \nabla_g \partial_{xx} u , \nabla_g u \rangle -\frac{1}{2}\Big ( \partial_x |\nabla_g u |^2 \varphi^2 \Big )|_{x=0}.
\nonumber\end{eqnarray}
Now, we use the boundary condition in~\eqref{pbpos}
to obtain that, on~$M$,
$$ f'(u)\nabla_g u=\nabla_g(f(u))=\nabla_g\partial_\nu u
=-\nabla_g \partial_x u.$$
Therefore,
\begin{equation}\label{AK}\begin{split}&
-\frac{1}{2}\int_{M} \Big ( \partial_x 
|\nabla_g u |^2 
\varphi^2 \Big )|_{x=0}
-\int_M \langle \nabla_g u_x ,
\nabla_g u \rangle \varphi^2\\ &\qquad\qquad=
\int_M f'(u)|\nabla_g u |^2 \varphi^2. \end{split}
\end{equation}
All in all, by collecting the results in~\eqref{AQ}, \eqref{AZ},
and~\eqref{AK}, we obtain that
\begin{equation}\label{AA}\begin{split}
&\int_{\mathcal M} \varphi^2   |\nabla |\nabla_g u | |^2 - 
\int_{\mathcal M} \varphi^2  \Big \{ \frac{1}{2} \partial_{xx} |\nabla_g u|^2 + |H_g u|^2 +Ric_g (\nabla_g u , \nabla_g u) \Big \} 
\\&\qquad+
\int_{\mathcal M} \varphi^2 \langle \nabla_g \partial_{xx} u , 
\nabla_g u \rangle + \int_{\mathcal{M}}|\nabla_g u|^2  |\nabla 
\varphi |^2  \geq 0. \end{split}\end{equation}
Also, we observe that
$$|\partial_x|\nabla_g u | |^2+
\langle  \nabla_g \partial_{xx} u , \nabla_g u \rangle-\frac{1}{2} \partial_{xx} |\nabla_g u|^2=$$
$$
|\partial_x|\nabla_g u | |^2-|\partial_x \nabla_g u |^2 \leq 0
$$
by the Cauchy-Schwarz inequality. 

Accordingly,
$$|\nabla|\nabla_g u||^2=
|\nabla_g|\nabla_g u||^2+|\partial_x|\nabla_g u||^2
\le
\frac12\partial_{xx}|\nabla_g u|^2-
\langle  \nabla_g \partial_{xx} u , \nabla_g u \rangle.$$
This and~\eqref{AA}
give \eqref{poincare}.
\end{proof}

\section{Proof of theorems 
\ref{liouville1}--\ref{flatness}}\label{s-2-2}

With~\eqref{poincare} at hand, 
one can prove theorems
\ref{liouville1}--\ref{flatness}. 

For this scope, first, we recall the following 
lemma, 
whose proof can be found in
section~2 of~\cite{FSV-man}. 

\begin{lemma}\label{signe}
For any smooth $\phi:\mathcal M \rightarrow\R$, we have that
\begin{equation}\label{POS}
|H_{\bar g}\phi|^2\ge\big|\nabla_ {\bar g}|\nabla_{\bar g} \phi|\big|^2\qquad
{\mbox{
almost everywhere.}}
\end{equation}
\end{lemma}

Moreover, we have the following result:

\begin{lemma}\label{lad}
Let~$u$ be a bounded solution of~\eqref{pbpos}.
Assume that
$$Ric_g \geq 0$$
and that~$Ric_g$ does not vanish identically
on~$M$.

Suppose that
\begin{equation}\label{91}{\mbox{
$Ric_g(\nabla_gu,\nabla_gu)$ vanishes
identically on~${\mathcal M}$.}}\end{equation}

Then,~$u$ is constant on~${\mathcal{M}}$.
\end{lemma}

\begin{proof}
By assumption,
we have that~$Ric_g$ is strictly positive definite
in a suitable non empty open set~$U\subseteq M$.

Then,~\eqref{91} gives that~$\nabla_g u$ vanishes
identically in~$U\times\R^+$.

This means that, for any fixed~$x\in\R^+$,
the map~$U \ni y\mapsto u(x,y)$ does not depend on~$y$.
Accordingly,
there exists a function~$\tilde u:\R^+\rightarrow\R$ such that~$
u(x,y)=\tilde u(x)$, for any~$y\in U$.

Thus, from~\eqref{pbpos}, 
$$ 0=\Delta_{\bar g} u=\tilde u_{xx} \qquad{\mbox{
in $U\times\R^+$}}$$
and so
there exist~$a$, $b\in\R$ for which
$$ u(x,y)=\tilde u(x)=a+bx
\qquad{\mbox{
for any~$x\in\R^+$ and any~$y\in U$.}}
$$
Since~$u$ is bounded, we have that~$b=0$,
so~$u$ is constant in~$U\times\R^+$.

By
the unique continuation
principle (see Theorem 1.8 of \cite{Kazdan}), we have
that~$u$ is constant on~$M\times\R^+$.
\end{proof}

\subsection{Proof  of theorem \ref{liouville1}}\label{s-2-1-2}

Points in~${\mathcal{M}}$ will be denoted here as~$(x,y)$, 
with~$x\in \R^+$ and~$y\in M$.

Take $\varphi $ 
in 
\eqref{poincare} to be the function 
$$\varphi(x,y)=\phi(\frac{x}{R})$$
where $R>0$ and $\phi$ is a smooth cut-off, that is $\phi=0$ on 
$|x|\geq 2$ 
and $\phi=1$ on $|x|\leq 1$. 

We remark that this is an admissible test function, since~$M$
is assumed to be compact in theorem \ref{liouville1}.
Moreover, we remark that
\begin{equation}\label{A5}
|\nabla\varphi(x,y)|\le \frac{
\| \phi\|_{C^1 (\R)}\;\chi_{(0,2R)}(x)
}{R}.\end{equation}
Also, since $u$ is bounded, elliptic 
regularity gives
that $\nabla u$ is bounded in $M \times \R^+$. 

Therefore, using~\eqref{poincare},
lemma \ref{signe} and~\eqref{A5}, we obtain
\begin{equation}
\int_{M \times \R^+} \Big \{ Ric_g(\nabla_g u, \nabla_g u)  \Big \} \varphi^2 \leq \frac{C}{R^2}\int_{M \times (0,2R)} dV_{\bar g} \leq \frac{C}{R}  
\end{equation}
for some constant $C>0$. 

Sending $R\rightarrow+\infty$ and using the fact that $Ric_g \geq 
0$, we conclude that~$Ric_g(\nabla_g u, \nabla_g u)$ vanishes
identically.

Thus, by lemma~\ref{lad},
we deduce that $u$ is constant. 

\subsection{Proof of theorem \ref{liouville1NL}}  

We put coordinates~$x\in\R^+$ and~$y\in M$ for
points in~${\mathcal{M}}= M\times\R^+$.

Given a smooth and bounded~$u_o:M\rightarrow \R$,
we can define the harmonic extension~$\EE u_o: M\times\R^+
\rightarrow \R$ as the unique bounded function
solving
\begin{equation}\label{61}
\left\{
\begin{matrix}
\Delta_{\bar g}(\EE u_o)=0 & {\mbox{ in }}M\times\R^+,\\
\EE u_o=u_o & {\mbox{ on }}M\times\{ 0\}
.
\end{matrix}
\right.
\end{equation}
See Section 2.4 of~\cite{CSM} for furter details.

Then, we define
\begin{equation}\label{62}
\LL u_o:= \partial_\nu 
(\EE u_o)\big|_{x=0}.\end{equation}
We claim that, for any point in~$M\rightarrow \R$,
\begin{equation}\label{67}
-\partial_x (\EE u_o)\,=\,
\EE (\LL u_o).
\end{equation}
Indeed, by differentiating the PDE in~\eqref{61},
$$ \Delta_{\bar g}  \partial_x (\EE u_o)=0.$$
On the other hand,
$$ -\partial_x (\EE u_o)(0,y)=\partial_\nu(\EE u_o)(0,y)=
\LL u_o,$$
thanks to~\eqref{62}.

Moreover,~$\partial_x (\EE u_o)$ is bounded
by elliptic estimates, since so is~$u_o$.

Consequently, $-\partial_x (\EE u_o)$ is a bounded solution 
of~\eqref{61}
with~$u_o$ replaced by~$\LL u_o$.

Thus, by the uniqueness of bounded solutions of~\eqref{61},
we obtain~\eqref{67}.

By exploiting~\eqref{62} and~\eqref{67}, we see that
\begin{equation}\label{71}
\begin{split}
& \LL^2 u_o = 
\partial_\nu
\big( \EE ( \LL u_o)\big) \big|_{x=0}=
-\partial_x
\big( \EE ( \LL u_o)\big) \big|_{x=0}\\
&\qquad\;=
-\partial_x
\big( 
-\partial_x (\EE u_o)
\big) \big|_{x=0}=
\partial_{xx} (\EE u_o)\big|_{x=0}.
\end{split}\end{equation}
On the other hand,
using the~PDE in~\eqref{61},
$$ 0=\Delta_{\bar g} (\EE u_o)=
\Delta_{g} (\EE u_o)+\partial_{xx} (\EE u_o),$$
so~\eqref{71} becomes
$$ \LL^2 u_o (y) = 
\partial_{xx} (\EE u_o) (0,y)= -\Delta_{g} (\EE u_o)(0,y)
=-\Delta_g u_o (y),$$
for any~$y\in M$, that is
\begin{equation}\label{72}
\LL=(-\Delta_{g})^{1/2}.
\end{equation}
With these observations in hand, we now take~$u$
as in the statement of
theorem \ref{liouville1NL}
and we define~$v:= \EE u$.

{F}rom~\eqref{62} and~\eqref{72},
$$ \partial_\nu v\big|_{x=0}= \partial_\nu (\EE u)
\big|_{x=0}=\LL u=(-\Delta_{g})^{1/2} u.$$

Consequently, 
recalling~\eqref{uhdsiuuiuu111},
we obtain that~$v$ is a bounded solution of~\eqref{pbpos}.

Furthermore, the function $v$ is stable, thanks to~\eqref{staze}.

Hence $v$ is constant by theorem \ref{liouville1},
and so we obtain the desired result for $u=v|_{x=0}$. 

\subsection{Proof of theorem \ref{liouville2}}  
Given~$p=(m,x)\in M\times\R^+$, we define~$d_g(m)$ to be
the geodesic distance of~$m$ in~$M$ (with respect to a fixed
point) and
$$ d(p):=\sqrt{d_g(m)^2 +x^2}.$$
Let also~$\hat B_R:=\{ p\in M\times\R^+{\mbox{ s.t. }}
d(p)<R\}$, for any~$R>0$.
Notice that~$|\nabla_g u|\in L^\infty (M\times\R^+)$,
by elliptic estimates, and that~$\hat B_R\subseteq B_R\times[0,R]$,
where~$B_R$ is the corresponding geodesic ball in~$M$.

As a consequence, by our assumption on the volume of $B_R$,
we obtain
\begin{equation*}\begin{split}&
\int_{\hat B_R} |\nabla_g u|^2 \,dV_{\bar g}\le
\| \nabla_g u\|^2_{L^\infty(M\times\R^+)}
\int_{B_R\times [0,R]}\,dV_{\bar g}\\ &\qquad\;=
R\,\| \nabla_g u\|^2_{L^\infty(M\times\R^+)}
\int_{B_R}\,dV_{g}
\le C R(R+1) \,\| \nabla_g u\|^2_{L^\infty(M\times\R^+)}.\end{split}
\end{equation*}
That is, by changing name of~$C$,
\begin{equation}
\label{EN}
\int_{\hat B_R} |\nabla_g u|^2 \,dV_{\bar g}\le CR^2\qquad\quad
{\mbox{
for any $R\ge1$.}}
\end{equation}
Also, since~$d_g$ is a distance function
on~$M$ (see pages~34
and~123 of~\cite{Pet}), we have that
\begin{equation}\label{09}
|\nabla d(p)|=\frac{
\big| \big(d_g(m)\nabla_g d_g(m),\,x \big)\big|
}{d(p)} \le1.
\end{equation}
Also,
given~$R\ge 1$, we define
$$ \phi_R(p):=\left\{
\begin{matrix}
1 & {\mbox{ if $d(p)\le \sqrt R$,}}\\
(\log \sqrt R)^{-1} \big(\log R-\log(d(p))\big)
& {\mbox{ if $d(p)\in(\sqrt R,R)$,}}
\\
0 & {\mbox{ if $d(p)\ge R$.}}
\end{matrix}
\right.$$
Notice that (up to a set of zero~$V_{\bar g}$-measure)
\begin{equation*}
|\nabla\phi_R(p)|\le
\frac{
\chi_{\hat B_{R}\setminus \hat B_{\sqrt R}} (p)}{\log \sqrt R\; 
d(p)},
\end{equation*}
due to~\eqref{09}.

As a consequence, 
\begin{eqnarray*}
&&
(\log \sqrt R)^2 \int_{M \times \R^+} |\nabla_g u |^2 | \nabla 
\phi_R
|^2\,dV_{\bar g}\le
\int_{\hat B_R
\setminus \hat B_{\sqrt R}
}\frac{
|\nabla_g u (p)|^2
}{d(p)^2}\,dV_{\bar g}(p)
\\
&&\qquad=
\int_{\hat B_R
\setminus \hat B_{\sqrt R}
}
|\nabla_g u (p)|^2
\Big( 
\frac{1}{R^2}+\int_{d(p)}^{R}\frac{2\,dt}{t^3}
\Big)
\,dV_{\bar g}(p)
\\ &&\qquad\le
\frac{1}{R^2} \int_{\hat B_R}
|\nabla_g u (p)|^2
\,dV_{\bar g}(p)+
\int_{\sqrt{R}}^R
\int_{\hat B_t}
\frac{2|\nabla_g u (p)|^2}{t^3}
\,dV_{\bar g}(p)\,dt.
\end{eqnarray*}
Therefore, by~\eqref{EN},
$$
(\log \sqrt R)^2 \int_{M \times \R^+} |\nabla_g u |^2 | \nabla 
\phi_R
|^2\,dV_{\bar g}\le
C\left( 1
+ \int_{\sqrt{R}}^R
\frac{2\,dt}{t}
\right)\le 3C\log R
.$$
Consequently,
from~\eqref{poincare},
\begin{eqnarray}\label{fin qui}
&& \int_{M \times \R^+} \Big \{ Ric_g (\nabla_g u, \nabla_g  u) 
+|H_g
u|^2 -|\nabla_g |\nabla_g u||^2 \Big \} \phi^2_R \,
\le\,\frac{12 C}{\log 
R}. \end{eqnarray}

{F}rom this and~\eqref{POS},
we conclude that
$$
\int_{M \times \R^+} Ric_g (\nabla_g u, \nabla_g
u)\phi_R^2 
\,\le\,
\frac{12 C}{\log
R}. $$
By sending~$R\rightarrow+\infty$, we obtain
that~$Ric_g (\nabla_g u, 
\nabla_g  u)$ vanishes identically.

Hence,~$u$ is constant, thanks to lemma~\ref{lad},
proving
theorem \ref{liouville2}.

\subsection{Proof of theorem \ref{flatness}} 
The proof of
theorem~\ref{liouville2} can be carried out in this case
too, up to formula~\eqref{fin qui}.

Then,~\eqref{fin qui} in this case gives that
$$
\int_{M \times \R^+} \Big \{ |H_g
u|^2 -|\nabla_g |\nabla_g u||^2 \Big \} \phi^2_R \,
\le\,\frac{12C}{\log
R}.$$
By sending~$R\rightarrow+\infty$, and by recalling~\eqref{POS},
we conclude that~$|H_g u|$ is identically equal 
to~$|\nabla_g|\nabla_g u||$ on~$\big(
M\times\{x\}\big)\cap 
\{\nabla_g u\ne0\}$, for any fixed~$x>0$.

Consequently,
by lemma~5 of~\cite{FSV-man},
we have that
for any $k=1,\dots,n$ there exist $\kappa^k:M\rightarrow
\R$
such that
\begin{equation*}
\nabla_g\big( \nabla_g u\big)^k(p)=
\kappa^k(p) \nabla_g u(p)\qquad{\mbox{
for any $p\in \big( M\times\{x\}\big)\cap \{ \nabla_gu\ne0\}$.}}
\end{equation*}
{F}rom this and~\cite{FSV-man} (see the computation
starting there on formula~(23)), 
one concludes that
every connected component of~$
\{ y\in M,\,\,\,u(x,y)=c\}$
is a geodesic.

\section{The case of the hyperbolic space }\label{negative}

We now come to problem \eqref{pbneg}. Notice that up to now we assumed for the manifold $\mathcal M$ to be positively curved. We deal here with special equations on negatively curved manifolds. As a consequence, the geometric  formula \eqref{poincare} is not useful since the Ricci tensor does not have the good sign
and so we need a different strategy to deal with
the hyperbolic case.

For this, we will make use 
here of the fact that the manifold $\mathbb 
H^{n+1}$ with the metric $\bar g
\frac{|dy|^2+\abs{dx}^2}{x^2}
$ is conformal to $\R^{n+1}_+$ 
with the flat metric, and, in fact, $({\mathbb{H}}^{n+1},  
\bar g)$ is the main example of 
conformally compact Einstein manifold, as we discuss
in section \ref{Ein}
here below.

\subsection{Motivations and scattering theory}\label{Ein}
In order to justify the study of problem \eqref{pbneg}, we describe the link between problem \eqref{pbneg} and fractional order conformally covariant operators. 

Let $M$ be a
compact manifold of dimension $n$. Given a metric $h$ 
on $M$, the conformal class $[h]$ of $h$ is defined as the set of metrics $\hat h$ that can be written as $\hat h=f h$ for a positive conformal factor $f$. 

Let $\mathcal M$ be a smooth manifold of dimension $n+1$ with
boundary $\partial\mathcal M=M$. 

A function $\rho$ is a 
\emph{defining 
function} of $\partial \mathcal M$ in $\mathcal M$ if
$$\rho>0 \mbox{ in } \mathcal M, \quad \rho=0 \mbox{ on }\partial \mathcal M,\quad d\rho\neq 0 \mbox{ on } \partial \mathcal M $$
We say that $g$ is a \emph{conformally compact} metric on $X$ with conformal infinity $(M,[h])$ if there exists a defining function $\rho$ such that the manifold $(\bar{\mathcal M},\bar g)$ is compact for $\bar g=\rho^2 g$, and $\bar g|_M\in [h]$. 

If, in addition $(\mathcal M^{n+1}, g)$ is a conformally compact manifold and $Ric_g = -ng$, then we call $(\mathcal M^{n+1}, g)$ a \emph{conformally compact Einstein} manifold.

\medskip

Given a conformally compact, asymptotically hyperbolic manifold 
$(\mathcal M^{n+1}, g)$ and
a representative $\hat g$ in $[\hat g]$
on the conformal infinity $M$, there is a uniquely defining
function $\rho$ such that, on $M \times (0,\epsilon)$ in $\mathcal M$, $g$ has the normal form
$g = \rho^{-2}(d\rho^2 + g_\rho)$ where $g_\rho$ is a one 
parameter family of metrics on $M$ (see~\cite{GZ}
for precise statements and further details).
\medskip

In this setting, the scattering matrix of $M$ is defined as follows. Consider the following eigenvalue problem in $(\mathcal M,g)$, with Dirichlet boundary condition,
\begin{equation}
\label{1818}\left\{\begin{matrix}
-\Delta_g u_s-s(n-s)u_s& =0\mbox{ in } \mathcal M \\
u_s& = f \mbox{ on } M 
\end{matrix}\right.
\end{equation}
for $s\in\mathbb C$ and $f$ defined on $M$. 

Problem~\eqref{1818}
is solvable unless $s(n-s)$ belongs to the spectrum of $-\Delta_g$.  

However,
$$\sigma(-\Delta_g) = \left[(n/2)^2,\infty\right) \cup \sigma_{pp}(\Delta_g)$$
where the pure point spectrum $\sigma_{pp}(\Delta_g)$ (i.e., the 
set of $L^2$ eigenvalues), is finite and it is contained in
$\left(0,(n/2)^2\right)$. 

Moreover, given any  $f$ on $M$, Graham-Zworski \cite{GZ} obtained a 
meromorphic family of solutions $u_s=\mathcal{P}(s)f$ such that,
if $s \not\in n/2 + \mathbb{N}$, then
$$\mathcal{P}(s)f = F \rho ^{n-s} + H\rho^s.$$
And if $s = n/2 + \gamma$, $\gamma \in \mathbb N$,
$$\mathcal P(s)f = F\rho ^{n/2-\gamma} + H \rho^{ n/2+\gamma} \log \rho$$
where $F, H \in \mathcal C^\infty(X)$, $F|_M = f$, and $F,H$ mod $O(\rho^n)$ are even in $\rho$. 

It is worth mentioning that
in the second case $H|_M$ is locally determined by $f$ and $\hat g$. However,
in the first case, $H|_M$ is globally determined by $f$ and $g$. We are interested in the study of these nonlocal operators.

We define the scattering operator as $S(s)f = H|_M$, which is a 
meromorphic family of pseudo-differential operators in $Re(s)>n/2$ with poles at
$s=n/2+\mathbb N$ of finite rank residues. The relation between $f$ and $S(s)f$ is like that of the Dirichlet to Neumann operator in standard harmonic analysis. Note that the principal
symbol is
$$\sigma\lp S(s)\rp= 2^{n-2s}\;\frac{\Gamma(n/2-s)}{\Gamma(s-n/2)}\;\sigma\lp(-\Delta_g)^{s-n/2}\rp$$

The operators obtained when $s=n/2+\gamma$, $\gamma\in\mathbb N$ 
have been well studied. Indeed, at those values of $s$ the 
scattering matrix $S(s)$ has a simple pole of finite rank and its 
residue can be computed explicitly, namely
$$Res_{s=n/2 +\gamma}S(s) = c_\gamma P_{\gamma},\quad  c_\gamma = (-1)^\gamma[2^{2\gamma}\gamma!(\gamma - 1)!]^{-1}$$
and $P_{\gamma}$ are the conformally invariant powers of the 
Laplacian constructed by~\cite{FG, GJMS}.

In particular, when $\gamma=1$ we have the conformal Laplacian,
$$P_1=-\Delta +\frac{n-2}{4(n-1)} R$$
and when $\gamma=2$, the Paneitz operator
$$P_2=\Delta^2 +\delta\lp a_n Rg+b_n Ric\rp d+\tfrac{n-4}{2}Q^n$$

We can similarly define the following fractional order operators on $M$ of order $\gamma\in(0,1)$ as
$$P_\gamma f := d_\gamma S(n/2+\gamma) f,\quad d_\gamma=2^{2\gamma}\frac{\Gamma(\gamma)}{\Gamma(-\gamma)}.$$
It is important to mention that these operators are conformally covariant. Indeed, for a change of metric $g_u=u^{\frac{4}{n-2\gamma}}g_0$, we have
$$P^{g_u}_\gamma f= u^{-\frac{n+2\gamma}{n-2\gamma}} 
P_\gamma^{g_0} \lp u f\rp.$$

The following result, which
can be found in \cite{CG}, establishes a link between 
scattering theory on $\mathcal M$ and a local problem in the half-space. We provide the proof for sake of completeness. 
\begin{lemma} \label{prop-extension}
Fix $0<\gamma<1$ and let $s=\frac{n}{2}+\gamma$. Assume that~$u$ 
is a smooth solution of 
\begin{equation}\label{u}
\left\{
\begin{matrix}
-\Delta_{\bar g} u-s(n-s)u=0 \,\,\, \mbox{ in }\,\,\, \mathbb 
H^{n+1},\\
\partial_\nu u=v \,\,\, \mbox{ on }\,\,\, \partial
\mathbb H^{n+1}.
\end{matrix}\right.
\end{equation}
for some smooth function $v$ defined on $\partial \mathbb H^{n+1}$. 

Then the function $U=x^{s-n}u$ solves 
\be\label{CS}
\left\{\begin{matrix}
\mbox{div}\,(x^{1-2\gamma} \nabla U)= 0& \quad \mbox{for 
}y\in\R^n,\; x\in (0,+\infty)\\
U(0,.)= u|_{x=0},& \quad \mbox{in  
} \R^n \\
-\lim_{x\to 0}x^{1-2\gamma}\partial_x U=C v 
\end{matrix}\right.
\ee
for some constant $C$. 
\end{lemma}
\begin{proof}
By the results in \cite{GZ}, one has the following representation of $u$ in $\mathbb H^{n+1}$ 
$$
u = x^{n-s} u|_{x=0} + x^s \partial_\nu u .
$$
{F}rom this we deduce that 
$$U=u|_{x=0} +x^{2s-n}\partial_\nu u.$$
Since $s=\frac{n}{2}+\gamma$, we have $2s-n=2\gamma >0$ and then
$$U|_{x=0}=u|_{x=0}$$
and 
$$-\lim_{x\to 0}x^{1-2\gamma}\partial_x U= C \partial_\nu u =C v. $$
We now prove that $U$ satisfies the desired equation. This only comes from the conformality of the metric on $\mathbb H^{n+1}$ to the flat one in the half-space. Indeed, the conformal Laplacian is given by 
$$L_{\bar g} =-\Delta_{\bar g} +\frac{n-1}{4 n} R_{\bar g}$$
where $R_{\bar g}$ is the scalar curvature of $\mathbb H^{n+1}$, which is equal to $-n(n+1).$ On the other hand, if we have $h=e^{2w} \bar g$ for some function $w$ (i.e. the metrics $h$ and $g$ are conformal) then the conformal law of $L_{\bar g}$ is given by
$$L_{h}\psi =e^{-\frac{n+3}{2}w}L_{\bar g}(e^{\frac{n-1}{2}w}\psi)$$
for any smooth $\psi$.
 
In our case, we have $h=|dx|^2+|dy|^2$ the flat metric on $\mathbb R^{n+1}_+$ and $e^{w}=x$. 
Thus, using the conformal law, we have 
$$-\Delta_{\bar g} \psi= -x^2 \Delta \psi + (n-1) \partial_x \psi$$ 
for any $\psi$ smoothly on $\mathbb R^{n+1}_+. $

Plugging $\psi=u$ and using  equation \eqref{u} leads 
$$s(n-s)u= -x^2 \Delta u +(n-1) x\partial_x u .$$

Finally, plugging $U= x^{s-n}u$ leads to the equation 
$$\Delta U +\frac{1-2\gamma}{x}\partial_x U=0,$$
which is equivalent to $\mbox{div}\,(x^{1-2\gamma} \nabla U)= 0$. 
\end{proof}

\subsection{Proof of theorem \ref{thhyp}}\label{989}

Let $u$ be a solution as requested
in Theorem \ref{thhyp}. By 
Lemma~\ref{prop-extension}, 
the function $U$ satisfies in a weak sense
\be
\left\{\begin{matrix}
\mbox{div}\,(x^{1-2\gamma} \nabla U)= 0& \quad \mbox{for 
}y\in\R^2,\; x\in (0,+\infty)\\
U(0,.)=u|_{x=0}, \\
-\lim_{x\to 0}x^{1-2\gamma}\partial_x U=f(U).
\end{matrix}\right.
\ee
Notice that
either~$\partial_{y_2} U>0$ or~$f' \leq 0$, thanks to~\eqref{001}
and~\eqref{002}. Furthermore, since $u$ is bounded, $U$ is bounded close to $x=0$. Additionally, we have 
$$U= x^{\gamma -1} u.$$
This gives that~$U$ 
is bounded inside~$\mathbb H^{n+1}_+$.
So, since~$U$ agrees with~$u$ on~$\partial\mathbb H^{n+1}_+$,
we obtain that
$U$ is bounded in all of $\mathbb H^{n+1}_+. $ 

Therefore, by theorem~3 in~\cite{SV1}, we have that~$U(x,y)=
U_o (x,\omega\cdot y)$, for 
suitable~$U_o:[0,+\infty) \times \R\rightarrow\R$
and~$\omega\in S^1$.
This gives directly the desired result. 

\bibliographystyle{alpha}
\bibliography{bibman}

\end{document}